\documentclass[a4paper,10pt]{article}

\def\C{\ensuremath{{\bf C}}}
\def\R{\ensuremath{{\bf R}}}
\def\Z{\ensuremath{{\bf Z}}}

\def\N{\ensuremath{{\bf N}}}
\def\T{\ensuremath{{\bf T}}}
\def\B{\ensuremath{\mathcal{B}}}

\def\h{\ensuremath{H}}
\def\hl{\ensuremath{\mathcal{H}}}

\def\uhl{\ensuremath{{\hl}}}
\def\supp#1{\ensuremath{\text{supp}(#1)}}

\newtheorem{Theorem}{Theorem}
\newtheorem{Proposition}{Proposition}
\newtheorem{Lemma}{Lemma}
\newtheorem{Corollary}{Corollary}
\newtheorem{Definition}{Definition}
\newtheorem{Notation}{Notation}
\newtheorem{Example}{Example}
\newtheorem{Remark}{Remark}
\def\BEq{\begin{equation}}
\def\EEq{\end{equation}}
\def\BRem{\begin{Remark}}
\def\ERem{\end{Remark}}
\def\BEx{\begin{Example}}
\def\EEx{\end{Example}}
\def\BThe{\begin{Theorem}}
\def\EThe{\end{Theorem}}
\def\BNot{\begin{Notation}}
\def\ENot{\end{Notation}}
\def\BDef{\begin{Definition}}
\def\EDef{\end{Definition}}
\def\BPro{\begin{Proposition}}
\def\EPro{\end{Proposition}}
\def\BLem{\begin{Lemma}}
\def\ELem{\end{Lemma}}
\def\BCor{\begin{Corollary}}
\def\ECor{\end{Corollary}}
\def\BProof{\noindent{Proof.~}}
\def\EProof{\\}

\usepackage[utf8x]{inputenc}
\usepackage{amsmath}

\title{A weak local irregularity property in $S^\nu$ spaces}
\author{ \small \sc Marianne Clausel\thanks{Universit\'e de Lyon, CNRS, INSA de Lyon, Institut Camille Jordan UMR 5208, B\^atiment L. de Vinci, 20 av. Albert Einstein, F-69621 Villeurbanne Cedex, France.} and Samuel Nicolay\thanks{Universit\'e de Li\`ege, Institut de Math\'ematique, Grande Traverse, 12, B\^atiment B37, B-4000 Li\`ege (Sart-Tilman), Belgium.} \thanks{Corresponding author. Email: S.Nicolay@ulg.ac.be. Phone: +32(0)43669433. Fax: +32(0)43669547.}}

\begin{document}

\maketitle

\begin{abstract}
Although it has been shown that, from the prevalence point of view, the elements of the $S^\nu$ spaces are almost surely multifractal, we show here that they also almost surely satisfy a weak uniform irregularity property.
\end{abstract}

\section{Introduction}
The uniform regularity defined from the H\"older spaces $C^\alpha(\R^d)$ is one of the most popular concepts for the uniform regularity. It has been introduced to study smoothness properties of functions such as the Weierstra{\ss} function (see \cite{har:16}). Indeed, many ``historical'' functions share the same property (see \cite{jaf-nic:09}): there exist $H\in (0,1)$ and a constant $C>0$ such that the function $f$ satisfies on some interval $I$,
\begin{equation}\label{eq:unif rir1}
\forall x,y\in I,\, |f(x)-f(y)| \le C |x-y|^H\;,
\end{equation}
and
\begin{equation}\label{eq:unif rir2}
\forall x,y\in I, \sup_{(u,v)\in [x,y]^2}|f(x)-f(y)| \ge \frac{1}{C} |x-y|^H\;.
\end{equation}
It has been shown in \cite{clanic:pre} that this behavior is the typical behavior of the functions belonging to $C^H(\R^d)$, in the sense of prevalence (the notion of prevalence will be defined in the sequel, see section~\ref{subsec:prev}). In other words, almost every function of $C^H(\R^d)$ satisfies Relations~(\ref{eq:unif rir1}) and~(\ref{eq:unif rir2}).

However, in many cases, the functions do not satisfy Relations~(\ref{eq:unif rir1}) and~(\ref{eq:unif rir2}). It is in particular the case of the so-called multifractal functions (see Definition~\ref{def:f multifrac}), originally introduced in the context of turbulence and now used in many fields of science (see e.g.\ \cite{abr:02,arn:02,hal:86,jaf:01,sre:91}). It can be shown that in several functional spaces, almost every function (in the sense of prevalence) is multifractal (see \cite{fra:06,aub:07}).

In this paper, we aim at investigating the typical irregularity properties of the elements of the $S^\nu$ spaces. These functional spaces, defined in \cite{jaf:snu,aub:06}, give rise to an efficient multifractal formalism (see section~\ref{sec:multifrac models}), i.e.\ an heuristic method to study the pointwise regularity of a function from a global point of view. It has been shown in \cite{aub:07} that almost every element of $S^\nu$ is multifractal. We show here that, although they are multifractal, the elements of $S^\nu$ also satisfy almost surely a weak uniform irregularity property. To this end, we introduce a concept of weak uniform irregularity, define the irregularity exponent $\uhl(x_0)$ of a function at a given point $x_0$ and prove that, in the sense of prevalence, there exists $\alpha$ (depending on $S^\nu$) such that almost every function of $S^\nu$ satisfy
\[
 \uhl(x)=\alpha,
\]
for any $x$.

The paper is organized as follows. In the next section, we recall some definitions about the multifractal analysis and the $S^\nu$ spaces. Next, we define the local irregularity exponent. To prove our main result, we will need wavelet criteria for the irregularity; these are stated in Section~\ref{sec:criteria}. The prevalent result is obtained in the last section.

\section{Multifractal models}\label{sec:multifrac models}
The multifractal analysis aims to study the smoothness of very irregular signals. For such a function, it is meaningless to try to characterize its regularity at a given point, since the pointwise regularity can abruptly change. One rather tries to determine the so-called spectrum of singularities, which gives the size of the set of points that share the same pointwise regularity. We first have to precise what is meant by ``pointwise regularity'' and ``size of set''. The $S^\nu$ spaces have been introduced in order to provide an efficient multifractal formalism, i.e.\ a method that allows the computation of the multifractal spectrum in many practical cases. Although this technique does not always lead to the right spectrum, it has been shown that the $S^\nu$-based multifractal formalism gives the right answer for almost every element of $S^\nu$. Here the term ``almost every'' has to be clearly stated, since one can not use the usual Lebesgue measure for the infinite dimensional settings. Finally, we introduce another index of regularity, called the local irregularity exponent.

\subsection{Pointwise regularity, Hausdorff dimension and multifractal spectrum}
The notion of the pointwise regularity which we will use here is based on the characterization of the H\"older spaces $C^\alpha(x_0)$.
\BDef
Let $x\in\R^d$ and $\alpha>0$. A locally bounded function $f:\R^d\to \C$ belongs to $C^\alpha(x_0)$ if there exist $C>0$ and a polynomial $P$ of degree strictly lower than $\alpha$ such that, in a neighborhood of $x_0$,
\[
 |f(x)-P(x-x_0)| \le C |x-x_0|^\alpha.
\]
The H\"older exponent of $f$ at $x_0$ is the quantity
\[
 \h(x_0)=\sup\{\alpha>0: f\in C^\alpha(x_0)\}.
\]
\EDef

The multifractal analysis provides a description of the collection of the H\"older exponents of a function by way of the multifractal spectrum (also called the singularity spectrum), which associates to a value $h$ the size of the set of points for which the H\"older exponent is $h$. By size, one usually means Hausdorff dimension; this notion of dimension is usually preferred because it relies on a measure. We will only give the necessary definitions; the interested reader is referred to e.g.\ \cite{fal:90,mat:95}.
\BDef
Let $S$ be a Borelian subset of $\R^d$ and $\Gamma_\epsilon(S)$ be the collection of all the countable $\epsilon$-coverings of $S$, i.e.\ the collection of all the countable coverings of $S$ by sets whose diameter is lower than $\epsilon$. The $\delta$-dimensional Hausdorff measure of $S$ is
\[
 m^\delta(S)=\lim_{\epsilon\to 0} \inf_{(S_j)_j \in \Gamma_\epsilon(S)} \sum_j |S_j|^\delta,
\]
where $|S_j|$ denotes the diameter of the set $S_j$.
\EDef
It is easy to check that if $m^{\delta_0}(S)$ is finite, then $m^\delta(S)=\infty$ if $\delta<\delta_0$ and $m^\delta(S)=0$ if $\delta>\delta_0$. We are then naturally led to the following definition.
\BDef
The Hausdorff dimension $\dim S$ of a non-empty Borelian subset $S$ of $\R^d$ is given by
\[
 \dim S=\inf\{\delta: m^\delta(S)=0\}.
\]
One sets $\dim \emptyset=-\infty$.
\EDef
If $S$ is empty or uncountable, one has $\dim S=\sup\{\delta: m^\delta(S)=\infty\}$.
We are now able to introduce the multifractal spectrum of a function.
\BDef
The multifractal spectrum of a locally bounded function $f:\R^d\to \C$ is the application
\[
 d : (0,\infty]\to [-\infty,d] \quad h\mapsto \dim\{x\in \R^d: \h(x)=h\}.
\]
\EDef
Finally, let us introduce some usual denominations.
\BDef\label{def:f multifrac}
A function is called multifractal if the associated multifractal spectrum takes more than one real value on a support with non-zero Hausdorff dimension. A function that is not multifractal is a monofractal function. If the spectrum only takes one value, the function is said to be monoH\"older.
\EDef

\subsection{Bases of wavelets and $S^\nu$ spaces}
One of the remarkable aspects of the wavelets is that they provide bases of the space $L^2(\R^d)$. Historically, the first orthonormal wavelet basis is the Haar basis, constructed long before the introduction of the term ``wavelet''. In \cite{jaf:04}, it is shown that, under some suitable conditions, the multifractal spectrum of a function can be estimated from its wavelet coefficients. The $S^\nu$ spaces were introduced in \cite{jaf:snu} to improve the classical wavelet-based multifractal formalism. The $S^\nu$ spaces are closely related to the Besov spaces.

Let us briefly recall some definitions and notations (for more precisions, see e.g.\ \cite{mey:90,dau:92,mal:98}). Under some general assumptions, there exist a function $\phi$ and $2^d-1$ functions $(\psi^{(i)})_{1\le i < 2^d}$ called wavelets such that
\[
 \{\phi(x-k): k\in \Z^d\} \cup \{\psi^{(i)}(2^j-k): 1\le i<2^d, k\in \Z^d, j\in \N \}
\]
form a basis of $L^2(\R^d)$. A function $f\in L^2(\R^d)$ can be decomposed as follows,
\[
 f= \sum_{k\in \Z^d} C_k \phi(\cdot -k) + \sum_{j=1}^\infty \sum_{k\in\Z^d} \sum_{i=1}^{2^d-1} c_{j,k}^{(i)} \psi^{(i)} (2^j \cdot -k),
\]
where
\[
 C_k= \int_{\R^d} f(x) \phi(x-k) \,dx
\]
and
\begin{equation}\label{eq:cjki}
 c_{j,k}^{(i)} =  2^{dj} \int_{\R^d} f(x) \psi^{(i)} (2^j x -k) \,dx.
\end{equation}
Let us remark that we do not choose the (usual) $L^2$ normalization for the wavelets, but rather an $L^\infty$ normalization, which is better fitted to the study of the H\"olderian regularity. Expressions such as (\ref{eq:cjki}) can make sense in more general settings (e.g.\ if $f$ is a distribution). Hereafter, we will assume that the wavelets belong to $C^\gamma(\R^d)$ with $\gamma\ge\alpha+1$, and that the functions $\{\partial^s \phi\}_{|s|\le\gamma}$, $\{\partial^s \psi\}_{|s|\le\gamma}$ have fast decay. Moreover, for the sake of simplicity, when dealing with the $S^\nu$ spaces, we will suppose that the application $f$ is defined on the torus $\T^d=\R^d/\Z^d$. Let
\[
 \Lambda=\left\{(i,j,k): 1\le i<2^d, j\in \N, k\in\{0,\ldots,2^j-1\}^d\right\}\;.
\]
If $(i,j,k)\in \Lambda$, the periodized
\[
 \psi^{(i)}_{j,k} = \sum_{l\in \Z^d} \psi \big( 2^j(\cdot -l) -k \big)
\]
form a basis of the one-periodic functions in $L^2([0,1]^d)$. We will denote $(c_{j,k}^{(i)})_{(i,j,k)\in\Lambda}$ or $(c_\lambda)_{\lambda\in\Lambda}$ the wavelet coefficients of a function in $L^2([0,1]^d)$.

Let us now introduce the $S^\nu$ spaces.
\BDef
For a sequence $c=(c_\lambda)_{\lambda\in\Lambda}$, $C>0$ and $\alpha\in\R$, we define
\[
 E_j(C,\alpha)[c]= \{(i,k): |c_{j,k}^{(i)}\ge C 2^{-j\alpha}\}.
\]
The wavelet profile $\nu_c$ of $c$ is defined as
\[
 \nu_c(\alpha)= \lim_{\epsilon\to 0^+} \limsup_{j\to\infty} \frac{\log_2 \#E_j(1,\alpha+\epsilon)[c]}{j}\;,
\]
with $\alpha\in\R$. If $c$ represents the wavelet coefficients of a function $f$, one sets $E_j(C,\alpha)[f]=E_j(C,\alpha)[c]$ and $\nu_f=\nu_c$.
\EDef
Clearly, $\nu_c$ is non-decreasing, right-continuous and non-negative (lower than $d$) when not equal to  $-\infty$. It gives, in some way, the asymptotic behavior of the number of coefficients of $c$ that have a given order of magnitude.
\BDef\label{def:snu}
Let $\nu:\R\to \{-\infty\} \cup [0,d]$ be a non-decreasing, right-continuous function such that there exists $\alpha_0$ for which $\nu(\alpha)=-\infty$ if $\alpha<\alpha_0$ and $\nu(\alpha)\in [0,d]$ otherwise. A sequence distribution $f$ belongs to $S^\nu$ if its wavelet coefficients satisfy $\nu_f(\alpha)\le \nu(\alpha)$ for any $\alpha\in\R$.
\EDef
These spaces are robust, in the sense that their definition does not depend on the choice of the wavelet basis (see \cite{jaf:snu}). Roughly speaking, a function belongs to $S^\nu$ if for each scale $j$ and every $\alpha\in \R$, there are about $2^{\nu(\alpha)j}$ coefficients $c_{j,k}^{(i)}$ that are larger than $2^{-\alpha j}$. The following definition is equivalent to the preceding one (see \cite{aub:06}).
\BDef
If $\nu$ is a function as given in Definition~\ref{def:snu}, $f$ belongs to $S^\nu$ if for any $\alpha\in\R$, $\epsilon>0$ and $C>0$, there exists $J$ such that $j\ge J$ implies
\[
 \#E_j(C,\alpha)[f] \le 2^{(\nu(\alpha)+\epsilon)j}.
\]
\EDef

A topological framework for the spaces $S^\nu$ will also be needed (further details can be found in \cite{aub:06,aub:07}). The ancillary spaces are useful to obtain a structure of complete metric space on $S^\nu$.
\BDef
Let $m,n\in\N$; $f$ belongs to $E_{m,n}$ if there exist $C>0$ such that
\begin{equation}\label{eq:anc sp}
 \# E_j(C,\alpha_m)[f]\le C 2^{\nu(\alpha_m)+\epsilon_n)j},
\end{equation}
for any $j$, where $(\alpha_m)_m$ is any dense sequence in $\R$ and $(\epsilon_n)_n$ is decreasing to zero.
\EDef
If for $f\in E_{m,n}$, the infimum of the constants $C$ satisfying the inequality (\ref{eq:anc sp}) is noted $d_{m,n}(f,0)$, the distance $d_{m,n}(f,g)=d_{m,n}(f-g,0)$ makes $(E_{m,n},d_{m,n})$ a metric space. It can be shown that $S^\nu=\cap_{m,n} E_{m,n}$ and that with the distance
\[
 d(f,g)= \sum_{m,n\ge 0} 2^{-m-n} \frac{d_{m,n}(f,g)}{1+d_{m,n}(f,g)},
\]
the space $(S^\nu,d)$ is a complete separable metric space. The distance $d$ may depend on the sequences chosen in (\ref{eq:anc sp}), but the induced topology does not. The Borel $\sigma$-algebra relative to this topology will be denoted $\B(S^\nu)$.

\subsection{Prevalence of multifractal functions in $S^\nu$}\label{subsec:prev}
In a finite dimensional space, a property ``holds almost everywhere'' if the set of points for which this property is not satisfied vanishes for the Lebesgue measure. The Lebesgue measure has a preponderant role, as there is no other $\sigma$-finite translation invariant measure. Unfortunately, there is such measure in the infinite dimensional Banach spaces. The notion of prevalence provides the analogue of ``Lebesgue measure zero'' in complete metric vector spaces.

In \cite{chris:72}, to recover the notion of ``almost every'' in infinite vector spaces, a well-known characterization of Lebesgue measure zero subsets of $\R^d$ is generalized: In $\R^d$, a Borel set $B$ has Lebesgue measure zero if and only if there exists a compactly supported probability measure $\mu$ such that $\mu(B+x)=0$ for any $x\in\R^d$. This characterization can be turned into a definition in more general settings and leads to the concept of Haar-null set, which provides the analogue of ``Lebesgue measure zero'' set for infinite dimensional spaces.
\BDef\label{def:harr-null}
Let $E$ be a complete metric vector space. A Borel set $B\subset E$ is Haar-null if there exists a Borel probability measure, strictly positive on some compact set $K\subset E$ such that $\mu(B+x)=0$ for any $x\in E$. A subset of $E$ is Haar-null if it is included in a Haar-null Borel subset of $E$. The complement of a Haar-null set is called a prevalent set.
\EDef

In \cite{aub:07}, the prevalent behavior of almost every function of $S^\nu$ has been studied and it has been proved that almost every function of $S^\nu$ is multifractal. More precisely, an upper bound of the spectrum $d$ of any $f\in S^\nu$ is first given: if $\alpha_\text{max}=\inf_{h\ge \alpha_0} \{h/\nu(h)\}$ and
\[
 d_\nu(h)=\left\{\begin{tabular}{ll}
                  $h \sup_{h'\in (0,h]} \{\nu(h')/h'\}$ & if $h\le \alpha_\text{max}$ \\
                  $1$ & otherwise
                 \end{tabular}\right.,
\]
one has, for all $f\in S^\nu$ and $h\in \R$, $d(h)\le d_\nu(h)$. Then it is shown that the sets
\[
 \{f\in S^\nu: \text{$d(h)=d_\nu(h)$ if $h\le \alpha_\text{max}$ and $d(h)=-\infty$ otherwise}\}
\]
and
\[
 \{f\in S^\nu: \text{the pointwise regularity is almost everywhere $\alpha_\text{max}$}\}
\]
are prevalent. The ``typical elements'' (in the sense of the prevalence) of $S^\nu$ are thus multifractal and do not satisfy the classical law of the iterated logarithm almost everywhere.

\subsection{Local regularity exponents}
We introduce here another notion of local regularity, which will be used to state that a weak irregularity condition is satisfied for almost every multifractal function of $S^\nu$. To this end, we first need to recall different concepts of global regularity, based on the global H\"older regularity.

If $\Omega$ is an open subset of $\R^d$, for any $h\in\R^d$, we will denote by $\Omega_h$ the set
\[
 \Omega_h= \{ x\in\R^d : [x,([\alpha]+1)h]\subset \Omega\},
\]
where the value $\alpha>0$ will be implied by the context and $[\alpha]$ denotes the greatest integer lower than $\alpha$. We will need the classical notion of finite difference.
\BDef
Let $x,h\in\R^d$ and $f:\R^d\to \C$; the first order difference of $f$ is $\Delta_h^1 f(x)= f(x+h)-f(x)$. For $n\ge 2$, the difference of order $n$ is defined by
\[
 \Delta_h^n f(x)= \Delta_h^{n-1} \Delta_h^1 f (x).
\]
\EDef
We can now introduce the H\"older spaces $C^\alpha(\Omega)$.
\BDef
Let $\Omega$ be an open subset of $\R^d$ and $\alpha>0$; a bounded function $f$ defined on a subset of $\R^d$ belongs to $C^\alpha(\Omega)$ if there exist $C,r_0>0$ such that for any $r\le r_0$,
\[
 \sup_{|h|\le r} \|\Delta_h^{[\alpha]+1} f (x) \|_{L^\infty(\Omega_h)} \le C r^\alpha.
\]
The function $f$ is said to be uniformly H\"olderian on $\Omega$ if for some $\alpha>0$, $f$ belongs to $C^\alpha(\Omega)$.
\EDef
%\BDef
%The uniform H\"older exponent of a uniformly H\"older function $f$ on an open set $\Omega$ is defined as
%\[
% \lhl(\Omega) = \sup\{ \alpha>0: f\in C^\alpha(\Omega)\}.
%\]
%\EDef

We will also need a notion of uniform irregularity.
\BDef
Let $\Omega$ be an open subset of $\R^d$, $\alpha>0$ and $\beta\in\R$; a bounded function $f$ defined on a subset of $\R^d$ belongs to $I^\alpha (\Omega)$ if there exist $C,r_0>0$ such that for any $r\le r_0$,
\[
 \sup_{|h|\le r} \|\Delta_h^{[\alpha]+1} f (x) \|_{L^\infty(\Omega_h)} \ge C r^\alpha.
\]
A function that belongs to $I^\alpha(\Omega)$ is said to be uniformly irregular with exponent $\alpha$.
\EDef
Let us remark that the statement $f\in I^\alpha(\Omega)$ is not equivalent to $f\notin C^\alpha(\Omega)$. In the latter case, for any $C>0$, there exists a sequence $(r_n)_n$ (depending on $C)$ decreasing to zero for which
\[
 \sup_{|h|\le r_n} \|\Delta_h^{[\alpha]+1} f\|_{L^\infty(\Omega_h)} \ge C r_n^\alpha.
\]
We are thus naturally led to the following definition.
\BDef
Let $\Omega$ be an open subset of $\R^d$, $\alpha>0$ and $\beta\in\R$; a bounded function belong to $C^\alpha_w (\Omega)$ if $f$ is defined on $\Omega$ and if for any $C>0$, there exists a decreasing sequence $(r_n)_n$ decreasing to zero such that
\[
 \sup_{|h|\le r_n} \|\Delta_h^{[\alpha]+1} f \|_{L^\infty(\Omega)} \le C r_n^\alpha,
\]
for any $n\in\N$. The set $C^\alpha_{\beta,w}(\Omega)$ is simply denoted $C^\alpha_w(\Omega)$. A function that belong to $C^\alpha_w(\Omega)$ is said to be weakly uniformly H\"older with exponent $\alpha$ on $\Omega$.
\EDef
\BDef
The upper uniform H\"older exponent of a weakly uniformly H\"older function $f$ on an open set $\Omega$ is defined as
\[
 \uhl(\Omega) = \sup\{ \alpha>0: f\in C^\alpha_w(\Omega)\}.
\]
\EDef

We now define the local irregular exponent, using the same approach as in \cite{lev:02}.
\BDef
A sequence $(\Omega_n)_n$ of open subsets of $\R^d$ is decreasing to $x_0\in\R^d$ if
\begin{itemize}
 \item $m<n$ implies $\Omega_n \subset \Omega_m$,
 \item $|\Omega_n|\to 0$ as $n\to \infty$,
 \item $\cap_n \Omega_n= \{x_0\}$.
\end{itemize}
\EDef
The following lemma is needed.
\BLem
If $(\Omega_n)_n$ and $(\Omega'_n)_n$ are two sequences of open sets that decrease to $x_0$, then
\[
 \sup_{n\in\N} \{\uhl(\Omega_n)\}
 = \sup_{n\in\N} \{\uhl(\Omega'_n)\}.
\]
\ELem
\BProof
Let us suppose that $\sup_{n\in\N} \{\uhl(\Omega_n)\} > \sup_{n\in\N} \{\uhl(\Omega'_n)\}$. There exists an index $n_1$ such that $\uhl(\Omega_{n_1})> \sup_{n\in\N} \{\uhl(\Omega'_n)\}$. Now let $r>0$ be such that $B(x_0,r)\subset \Omega_{n_1}$; since $(\Omega_n')_n$ is decreasing to $x_0$, there exists an index $n_2$ such that $\Omega'_{n_2}\subset B(x_0,r)$. One thus have
\[
 \uhl(\Omega'_{n_2}) \ge \uhl(\Omega_{n_1}) > \sup_{n\in\N} \{\uhl(\Omega'_n)\},
\]
which leads to a contradiction.
\EProof
\BDef
If $f$ is a bounded function, the local irregularity exponent of $f$ at $x_0$ is
\[
 \uhl(x_0)= \sup_n\{ \uhl(\Omega_n)\},
\]
where $(\Omega_n)_n$ is a sequence of open sets decreasing to $x_0$.
\EDef

\section{Wavelet criteria for the uniform irregularity and characterization of the local irregularity exponent}\label{sec:criteria}
In this section, we first give necessary and sufficient conditions for a function to belong to $I^\alpha(\Omega)$. Next we characterize the local irregularity exponent in terms of wavelet coefficients.

\subsection{Wavelet criteria for the uniform irregularity}
In what follows, we will assume that the multiresolution analysis is compactly supported (see \cite{dau:88}). The following result is shown in \cite{jaf:04}: in $\R$, if the wavelet basis belongs to $C^M(\R)$, there exist a fast decaying function $\Psi_M$ such that $\psi=\Delta_{1/2}^M \Psi_M$. Furthermore, the function $\Psi_M$ can be picked compactly supported with support included in this of $\psi$. In $\R^d$, we will use the tensor product wavelet basis (see \cite{mey:90,dau:92}),
\[
 \psi^{(i)}(x)=\Psi^{(1)}(x_1) \cdots \Psi^{(d)}(x_d),
\]
where $\Psi^{(i)}$ ($i\in \{1,\ldots,d\}$) are either $\psi$ or $\phi$, but at least one of them must equal $\psi$. We will also use the following notations: given $j\in\N$, $\Omega$ an open subset of $\R^d$ and a family of wavelets $\psi^{(i)}_{j,k}$, we set
\[
 I_j=\{ (i,k): \supp{\psi^{(i)}_{j,k}} \subset \Omega\}
\]
and
\[
 \| c^{(\cdot)}_{j,\cdot} \|_{\ell^\infty(\Omega)}= \sup_{I_j} |c^{(i)}_{j,k}|.
\]

The uniform regularity of a function is related to the decay rate of its wavelet coefficients (see \cite{mey:90}). Let $f$ be a bounded function and $\alpha\in (0,1)$; $f$ belongs to $\dot{C}^\alpha(\Omega)$ (where $\dot{C}^\alpha(\Omega)$ denotes the homogeneous version of the H\"older space $C^\alpha(\Omega)$) if and only if there exists $C>0$ such that for any $(j,k)$ such that for any $j\ge 0$,
\begin{equation}\label{eq:dotC}
 \|c^{(\cdot)}_{j,\cdot}\|_{\ell^\infty(\Omega)} \le C 2^{-\alpha j}.
\end{equation}
The following result gives a sufficient condition to belong to $I^\alpha(\Omega)$.
\BThe\label{thm:in I}
Let $\alpha>0$, $f\in C^\alpha(\Omega)$ and set $M=[\alpha]+1$. If there exists $C>0$ and $\gamma>1$ such that, for any $j\ge 0$,
\begin{equation}\label{eq:hyp cn}
 \max \{
 \sup_{j\le l \le j+ \log_2 j} \| c^{(\cdot)}_{l,\cdot} \|_{\ell^\infty(\Omega)} ,
 2^{-jM} \!\!\! \sup_{j-\log_2 j \le l \le j} (2^{lM} \| c^{(\cdot)}_{l,\cdot} \|_{\ell^\infty(\Omega)}) \}
 \ge C 2^{-j\alpha} j^\gamma,
\end{equation}
then $f\in I^\alpha(\Omega)$.

Now, if $f$ belongs to $I^\alpha(\Omega)$, there exist $C>0$ and $\beta\in (0,1)$ such that for any integer $j\ge 0$,
\begin{equation}\label{eq:hyp cs}
 \max \{
 \sup_{j\le l \le j+ \log_2 j} \| c^{(\cdot)}_{l,\cdot} \|_{\ell^\infty(\Omega)} ,
 2^{-jM} j^\beta \!\!\! \sup_{j-\log_2 j \le l \le j} (2^{lM} \| c^{(\cdot)}_{l,\cdot} \|_{\ell^\infty(\Omega)}) \}
 \ge C 2^{-j\alpha}.
\end{equation}
\EThe
\BProof
To prove the first part of the theorem, let us suppose that $f\in C^\alpha_w(\R^d)$ and let $C>0$. As shown in \cite{clanic:pre}, there exists some increasing sequence of integers $(j_n)_n$ such that for any $n\in\N$ and any $j\ge j_n$,
\begin{equation}\label{eq:cn}
 \sup_{|h|\le 2^{-j}} \|\Delta_h^M f \|_{L^\infty(\Omega_h)} \le C 2^{-j_n \alpha}.
\end{equation}
 Let us show that this inequality leads to a contradiction.

For the sake of simplicity, let us suppose that $\Psi^{(1)}=\psi$ (let us recall that one of the $\Psi^{(i)}$ is $\psi$). We then have, by definition of the wavelet coefficients, for any $j\ge 0$ and any $(i,k)\in I_j$,
\begin{eqnarray*}
 c^{(i)}_{j,k} &=& 2^{jd} \int_{\R^d} f(x) \Psi^{(1)} (2^j x_1-k_1) \cdots \Psi^{(d)} (2^j x_d-k_d) \,dx \\
               &=& 2^{jd} \int_{\R^d} f(x) \Delta_{1/2}^M \Psi_M (2^j x_1-k_1) \cdots  \Psi^{(d)} (2^j x_d-k_d) \, dx \\
               &=& 2^{jd} \int_{\R^d} \Delta_{1/2^{j+1} e_1}^M f (x) \Psi_M (2^j x_1-k_1) \cdots \Psi^{(d)} (2^j x_d-k_d) \, dx,
\end{eqnarray*}
where $e_1=(1,0,\ldots,0)$. Using the assumptions on the support of $\Psi_M$ and the definition of $I_j$, one has, for any $n\in\N$, any $j\ge j_n$ and any $(i,k)\in I_j$,
\begin{eqnarray*}
% |c^{(i)}_{j,k}| &\le& 2^{jd} \int_{\R^d} |\Delta_{1/2^{j+1}e_1}^M f (x)| |\Psi_M (2^j x_1-k_1) \cdots \Psi^{(d)} (2^j x_d-k_d)| \, dx \\
 |c^{(i)}_{j,k}| &\le& 2^{jd} \int_\Omega |\Delta_{1/2^{j+1}e_1}^M f (x)| |\Psi_M (2^j x_1-k_1) \cdots \Psi^{(d)} (2^j x_d-k_d)| \, dx \\
                 &\le& C 2^{jd} 2^{-j_n\alpha} \int_{\R^d} |\Psi_M (2^j x_1-k_1) \cdots \Psi^{(d)} (2^j x_d-k_d)| \, dx \\
                 &=& C 2^{-j_n\alpha} \| \Psi_M \otimes \cdots \otimes \Psi^{(d)} \|_{L^1(\R^d)},
\end{eqnarray*}
by using Relation (\ref{eq:cn}).

For $n\in\N$, let $l_n=j_n+ \gamma [\log_2 j_n]$; for $n$ sufficiently large, one has $l_n- \log_2 l_n\ge j_n$. Therefore, the following relations hold for $n$ sufficiently large (and any $(i,k)\in I_j$),
\[
 \sup_{l_n \le j \le l_n+\log_2 l_n} |c^{(i)}_{j,k}| \le C 2^{-j_n\alpha} \le C 2^{-l_n\alpha} l_n^{\gamma'}
\]
and
\[
 \sup_{l_n-\log_2 l_n \le j \le l_n} 2^{jM} |c^{(i)}_{j,k}| \le C 2^{l_n M} 2^{-j_n \alpha} \le C 2^{l_n M} 2^{-l_n \alpha} l_n^{\gamma'},
\]
which is in contradiction with the relation (\ref{eq:hyp cn}).

To prove the second part of the theorem we will use the following result (see \cite{mey:90}): Let $f\in C^\gamma(\Omega)$; since $C^\gamma(\Omega)\subset B^0_{\infty,\infty}(\Omega) \cap \dot{C}^\gamma(\Omega)$, the wavelet characterizations of these two functional spaces lead to the existence of a constant $C>0$ that does not depend on the function such that for any $h\in\R^d$ and any $x\in\Omega_h$, one has
\begin{equation}\label{eq:meya}
 |\Delta_h^{[\gamma]+1} f(x)| \le C \sup_{j\in \N} \| c^{(\cdot)}_{j,\cdot} \|_{l^\infty(\Omega)}\;,
\end{equation}
and
\begin{equation}\label{eq:meyb}
 |\Delta_h^{[\gamma]+1} f(x)| \le C |h|^\gamma \sup_{j\in \N} \{ 2^{j\gamma} \| c^{(\cdot)}_{j,\cdot} \|_{l^\infty(\Omega)}\}.
\end{equation}
Let us assume $f\in I^\alpha(\Omega)$ and that Property (\ref{eq:hyp cs}) is not satisfied. In this case, for any $C>0$ and any $\beta\in (0,1)$, there exists an increasing sequence of integers $(j_n)_n$ such that, for any $n\in\N$,
\begin{equation}\label{eq:hyp cs contra}
 \begin{array}{cc}
 \displaystyle \max \{
 \sup_{j_n\le l \le j_n+ \log_2 j_n} \| c^{(\cdot)}_{l,\cdot} \|_{\ell^\infty(\Omega)} ,\\
 & \displaystyle \hspace{-100pt}
 2^{-j_n M}j^\beta \!\!\! \sup_{j_n-\log_2 j_n \le l \le j_n} (2^{lM} \| c^{(\cdot)}_{l,\cdot} \|_{\ell^\infty(\Omega)}) \}
 \end{array}
 \le C 2^{-j_n \alpha}.
\end{equation}
Let us fix $C>0$, $x\in\Omega_h$ and let $n_0\in\N$, $h\in\R^d$ be such that $|h|\le 2^{-j n_0}$. By definition of $I^\alpha(\Omega)$, we have to show that $f\in C^\alpha_w(\Omega)$. We will use the following notations:
\[
 f_{-1}= \sum_{k\in \Z^d} C_k \phi(\cdot -k), \qquad f_j= \sum_{i=1}^{2^d-1} \sum_{k\in \Z^d} c^{(i)}_{j,k} \psi(2^j \cdot -k),
\]
with $j\ge 0$. Since $f$ is uniformly H\"older, $f_j$ and $\sum_{j\ge -1} f_j$ converge uniformly on any compact set and
\[
 \Delta_h^M f= \sum_{j\ge -1} \Delta_h^M f_j.
\]

We first consider the function $g_1=\sum_{j=-1}^{j_{n_0}} f_j$. Let us fix $\gamma\in ([\alpha]+1-\beta,[\alpha]+1)$. The regularity of the wavelets and property~(\ref{eq:meyb}) imply the existence of $C>0$ not depending on $n_0$, $x$ and $h$ such that
\begin{equation}\label{eq:g1}
 |\Delta_h^M g_1(x)| \le C |h|^\gamma \sup_{l\le j_{n_0}} \left(2^{l\gamma} \| c^{(\cdot)}_{l,\cdot} \|_{\ell^\infty(\Omega)}\right).
\end{equation}
Since $f\in C^\alpha(\Omega)$, there exists a constant such that, for any $n\in\N$,
\[
 \sup_{l\le j_{n_0}- \log_2 j_{n_0}} (2^{l \gamma} \|c^{(\cdot)}_{l,\cdot} \|_{\ell^\infty(\Omega)})
 \le C 2^{(j_{n_0}-\log_2 j_{n_0}) (\gamma-\alpha)}
 = C' \frac{2^{j_{n_0} (\gamma-\alpha)}} {j_{n_0}^{\gamma-\alpha}}.
\]
On the other hand, we have, using relation~(\ref{eq:hyp cs contra}),
\begin{eqnarray*}
 \lefteqn{\sup_{j_{n_0}-\log_2 j_{n_0} \le l \le j_{n_0}} (2^{l \gamma} \| c^{(\cdot)}_{l,\cdot} \|_{\ell^\infty(\Omega)})} \\
 &=& \sup_{j_{n_0}-\log_2 j_{n_0} \le l \le j_{n_0}} (2^{l (\gamma-M)} 2^{lM} \| c^{(\cdot)}_{l,\cdot} \|_{\ell^\infty(\Omega)}) \\
 &\le& C 2^{(j_{n_0} -\log_2 j_{n_0}) (\gamma-M)} 2^{j_{n_0} (M-\alpha)}/ j_{n_0}^{\beta} \\
 &=& C' j_{n_0}^{M-\gamma} 2^{j_{n_0} (\gamma-\alpha)}/j_{n_0}^\beta \\
 &\le& C' 2^{j_{n_0} (\gamma-\alpha)} .
\end{eqnarray*}
This implies that for $n_0$ sufficiently large (since $0<M-\gamma\le \beta$),
\[
 \sup_{l\le j_{n_0}} (2^{l \gamma} \|c^{(\cdot)}_{l,\cdot} \|_{\ell^\infty(\Omega)})
 \le C 2^{j_{n_0} (\gamma-\alpha)},
\]
and hence, using inequality~(\ref{eq:g1}),
\begin{equation}\label{eq:g1b}
 |\Delta_h^M g_1 (x)| \le C 2^{-j_{n_0} \alpha},
\end{equation}
since $h$ has been chosen adequately.

Let us now consider $g_2= \sum_{j> j_{n_0}} f_j$. Property~(\ref{eq:meya}) applied to $g_2$ directly gives the following relation,
\begin{equation}\label{eq:g2}
 |\Delta_h^M g_2| \le C \sup_{l>j_{n_0}} (\|c^{(\cdot)}_{l,\cdot} \|_{l^\infty(\Omega)}).
\end{equation}
Once again, since $f\in C^\alpha(\Omega)$, we have
\[
 \sup_{l\ge j_{n_0} +\log_2 j_{n_0}} \| c^{(\cdot)}_{l,\cdot}\|_{\ell^\infty(\Omega)} \le C 2^{-(j_{n_0}+\log_2 j_{n_0})\alpha} = C' \frac{2^{-j_{n_0} \alpha}}{j_{n_0}^\alpha}.
\]
Now, Relation~(\ref{eq:hyp cs contra}) implies
\[
 \sup_{j_{n_0} \le l \le j_{n_0} +\log_2 j_{n_0}} \|c^{(\cdot)}_{l,\cdot} \|_{\ell^\infty(\Omega)} \le C 2^{-j_{n_0} \alpha}.
\]
These last inequalities lead to the following relation for $n_0$ sufficiently large,
\begin{equation}\label{eq:g2b}
 |\Delta_h^M g_2 (x)| \le C 2^{-j_{n_0} \alpha},
\end{equation}
thanks to relation~(\ref{eq:g2}).

Putting relations~(\ref{eq:g1b}) and (\ref{eq:g2b}) together, we obtain
\[
 |\Delta_h^M f (x)| = |\Delta_h^M (g_1+g_2) (x)| \le C 2^{-j_{n_0} \alpha}
\]
for $n_0$ sufficiently large, that is $f\in C^\alpha_w (\Omega)$, which is impossible, since $f\in I^\alpha(\Omega)$.
\EProof

\subsection{Characterization of the local irregularity exponent}
The preceding result leads to the following characterization of the irregularity exponent.
\BCor\label{cor:carac}
Let $\alpha>0$; if $f\in C^\alpha(\Omega)$, then the irregularity exponent of $f$ on $\Omega$ equals $\alpha$ if and only if
\[
 \lim_{j\to \infty} \frac{\displaystyle
 \log_2 \max\{ \sup_{j\le l\le j+\log_2 j} \| c^{(\cdot)}_{l,\cdot} \|_{\ell^\infty(\Omega)} ,
 2^{-jM} \sup_{j-\log_2 j\le l \le j} (2^{lM} \| c^{(\cdot)}_{l,\cdot} \|_{\ell^\infty(\Omega)})
 \}}{-j} =\alpha,
\]
where $M=[\alpha]+1$.
\ECor
\BProof
Theorem~\ref{thm:in I} directly yields that if $f\in C^\alpha(\Omega)$, the irregularity exponent of $f$ on $\Omega$ equals $\alpha$ if and only if
\[
 \limsup_{j\to \infty} \frac{\displaystyle
 \log_2 \max\{ \sup_{j\le l\le j+\log_2 j} \| c^{(\cdot)}_{l,\cdot} \|_{\ell^\infty(\Omega)} ,
 2^{-jM}\!\!\! \sup_{j-\log_2 j\le l \le j} (2^{lM} \| c^{(\cdot)}_{l,\cdot} \|_{\ell^\infty(\Omega)})
 \}}{-j} =\alpha.
\]
We have to prove that the $\limsup$ can be replaced by a limit. Since $f\in C^\alpha(\Omega)$, the relation~(\ref{eq:dotC}) implies that
\[
 \liminf_{j\to \infty} \frac{\displaystyle
 \log_2 \max\{ \sup_{j\le l\le j+\log_2 j} \| c^{(\cdot)}_{l,\cdot} \|_{\ell^\infty(\Omega)} ,
 2^{-jM}\!\!\! \sup_{j-\log_2 j\le l \le j} (2^{lM} \| c^{(\cdot)}_{l,\cdot} \|_{\ell^\infty(\Omega)})
 \}}{-j} =\alpha.
\]
is always larger than
\[
 \liminf_{j\to\infty} \frac{\log_2 \max \{
 \sup_{j\le l\le j+\log_2 j} 2^{-l\alpha}, 2^{-jM} \sup_{j-\log_2 j\le l \le j} 2^{l(M-\alpha)}
 \}}{-j} =\alpha,
\]
which is sufficient to conclude.
\EProof
We can now state the local version of the previous result.
\BThe\label{thm:car H}
Let $\alpha>0$; if $f\in C^\alpha(\Omega)$ and $x_0\in \R^d$, then
\[
 \uhl(x_0)= \alpha
\]
if and only if
\[
 \lim_{r\to 0} \lim_{j\to \infty} \frac{
 \begin{array}{cc}
 \displaystyle
 \log_2 \max\{ \sup_{j\le l\le j+\log_2 j} \| c^{(\cdot)}_{l,\cdot} \|_{\ell^\infty(B(x_0,r))} ,\\
 & \displaystyle \hspace{-100pt}
 2^{-jM}\!\!\! \sup_{j-\log_2 j\le l \le j} (2^{lM} \| c^{(\cdot)}_{l,\cdot} \|_{\ell^\infty(B(x_0,r))})
 \}
 \end{array}}{-j} =\alpha,
\]
where $M=[\alpha]+1$.
\EThe

\section{A prevalent result on the $S^\nu$ spaces}
Prevalence supplies a natural definition of ``almost every'' which is translation invariant and where no specific measure plays a particular role. To prove our prevalent result, we first need to introduce the method we will use. We will also need some properties obtained in \cite{aub:07}.
\subsection{The stochastic process technique}
Our result concerning the prevalence relies on the stochastic process technique. Let us recall that a random element $X$ on a complete metric space $E$ is a measurable mapping $X$ defined on a probability space $(\Omega,\mathcal{A},P)$ with values in $E$. For a random element on $E$, one can define a probability on $E$ by the formula
\[
 P_X(A)=P\{ X\in A\}.
\]
Replacing the measure $\mu$ in Definition~\ref{def:harr-null} of a Haar-null set with $\mu=P_X$, we see that in order to prove that a set is Haar-null, it is sufficient to check that for any $f\in E$,
\[
 P_X(A+f)=0.
\]

The stochastic process that we will use here is a random wavelet series associated in a proper to $\nu$. To this end, for each $j\ge 0$, let us define as in \cite{aub:07},
\begin{equation}\label{eq:F_j}
 F_j(\alpha)= \left\{\begin{tabular}{ll}
                      $0$ & if $\alpha<\alpha_\text{min}$ ,\\
                      $2^{-jd} \sup\{j^2, 2^{j\nu(\alpha)}\}$ & if $\alpha \ge \alpha_\text{min}$ .
                     \end{tabular}\right.
\end{equation}
Since $F_j$ is non-decreasing and piecewise continuous, it is the repartition function of some probability distribution associated to a probability law $\rho_j$ supported on $[\alpha_\text{min},\infty]$, where
\[
 \alpha_\text{min}= \inf \{\alpha: \nu(\alpha)\ge 0\}.
\]
The following remark is made in \cite{aub:07}. If $\rho_j$ is the probability distribution whose repartition function $F_j$ is defined by~(\ref{eq:F_j}), then there exists some sequence of random numbers $(c_\lambda)_{\lambda\in\Lambda}$ with independent phase and moduli such that for any $j$, $\rho_j$ is the common law of $-\log_2 |c^{(i)}_{j,k}|/j$ and satisfies the two following conditions:
\begin{equation}\label{eq:soph:cond1}
 \forall \alpha\in\R,\qquad \lim_{j\to\infty} \frac{2^{jd} \rho_j((-\infty,\alpha])}{j} = \nu(\alpha),
\end{equation}
and
\begin{equation}\label{eq:soph:cond2}
 \forall \alpha\ge \alpha_\text{min},\qquad 2^{jd} \rho_j((-\infty,\alpha]) \ge j^2.
\end{equation}
Starting from these results, we will use the following random wavelet series associated to $\nu$,
\begin{equation}\label{eq:rws}
 X_\nu= \sum_{\lambda\in\Lambda} c_\lambda \psi_\lambda.
\end{equation}
It is shown in \cite{aub:06} that the metric topology $d(\cdot,\cdot)$ on $S^\nu$ makes it a Polish space, which is a very good framework for prevalence. Moreover, as proved in \cite{aub:07}, the measure $P_X$ is a Borel measure (relatively to this topology). In other words, we can use the stochastic process technique with $X_\nu$.

\subsection{Prevalent irregularity properties in $S^\nu$}
We are now ready to prove the following result.
\BThe
The following set is prevalent in $S^\nu$,
\[
 \{f: \text{$\forall x\in\T^d$, the local irregularity exponent of $f$ at $x$ is $\uhl(x)=\alpha_\text{min}$}\},
\]
where
\[
 \alpha_\text{min}= \inf \{\alpha: \nu(\alpha)\ge 0\}.
\]
\EThe
From this point of view, the ``typical elements'' of $S^\nu$ are multifractal and satisfy a weak uniform irregularity property.

From what precede, it is sufficient to show the following result.
\BPro
Let $X_\nu$ be the random wavelet series defined by~(\ref{eq:rws}). Then, for any $f\in S^\nu$, the local irregularity exponent $\uhl(x)$ of $f+X_\nu$ at $x$ is equal to $\alpha_\text{min}$ almost surely.
\EPro
\BProof
Since $S^\nu \subset C^{\alpha_\text{min}}(\T^d)$ then for any $f\in S^\nu$, $f+X_\nu \in C^{\alpha_\text{min}}(\T^d)$ almost surely. Let us now fix $m\in\N$ and define, for any $r\in \Z^d$,
\[
 T_{r,m} =  \prod_{n=1}^d ( \frac{r_n}{2^m}, \frac{r_n +1}{2^m})
\]
so that $\T^d= \cup_{r} \overline{T_{r,m}}$. We aim at showing that the equality
\[
 \uhl(T_{r,m})=\alpha_\text{min}
\]
holds almost surely for any $r,m$; in this case, Theorem~\ref{thm:car H} directly yields the required result.

For $\lambda\in \Lambda$, we will denote as usual $c_\lambda$ the wavelet coefficients associated to $X_\nu$ and $d_\lambda$ the wavelet coefficients associated to $f\in S^\nu$. Let us first remark that if for some fixed $\lambda\in \Lambda$, $\Re (c_\lambda \overline{d_\lambda})\le 0$ i.e.\ $c_\lambda$ is in the complex half-plane opposite to $d_\lambda$, then
\[
 |c_\lambda-d_\lambda| \ge |c_\lambda|
\]
Therefore, for any $\lambda\in \Lambda$,
\begin{equation}\label{eq:proba 1/2}
 P(f:|c_\lambda -d_\lambda|\ge |c_\lambda|) \ge P (f:\Re(c_\lambda \overline{d_\lambda})\le 0) \ge 1/2.
\end{equation}
If we define for any $N\in\N$ and $\lambda=(i,j,k)\in \Lambda$,
\[
 A_{N,\lambda} =\{f: (\exists k'\in k+[0,N]^d: |c^{(i)}_{j,k'}-d^{(i)}_{j,k'}|\ge |c^{(i)}_{j,k'}|)\},
\]
then, thanks to the independence of the wavelet coefficients of $X_\nu$ and inequality~(\ref{eq:proba 1/2}), we have
\[
 P(A_{N,\lambda}) \ge 1- 2^{-Nd}
\]
for any $N$ and $\lambda$.

Now, for $N\in\N$, let us set
\[
 B_{r,m,n} = \{f: \| c^{(\cdot)}_{j,\cdot} \|_{\ell^\infty(T_{r,m})} \ge 2^{-j (\alpha_\text{min} +1/n)} \}.
\]
If for any $i$, we have $\supp{\psi^{(i)}}\subset [0,M]^d$, then $\supp{\psi_\lambda}\subset T_{r,m}$ if and only if, for any $l\in\{1,\ldots,d\}$, we have $r_l 2^{j-m} \le k_l \le (r_l+1) 2^{j-m}-M$. Therefore, for any fixed $j\ge m$,
\[
 P(B_{r,m,n})=1- (1- 2^{j(\nu(\alpha_\text{min} +1/n)-d)})^{2^{d(j-m)}}
 \ge 1-\exp(-2^{-md}2^{j\nu(\alpha_\text{min}+1/n)}).
\]
Let us choose $N=2^{j-m}$; for any $m,n$ and any $j\ge m$, we have
\begin{eqnarray*}
 \lefteqn{P(f: \bigcap_r \| c^{(\cdot)}_{j,\cdot} -d^{(\cdot)}_{j,\cdot} \|_{\ell^\infty(T_{r,m})} \ge 2^{-j (\alpha_\text{min}+1/n)})} \\
 &\ge& \big(P(f: (\bigcup_{(i,k)\in I_j} A_{2^{j-m},\lambda}) \cap B_{2r,m+1,n})\big)^{2^m} \\
 &\ge& \big(1-2^{-2^{d(j-m)}} -\exp(-2^{-md} 2^{j\nu(\alpha_\text{min}+1/n)})\big)^{2^{md}}.
\end{eqnarray*}
Moreover, for any $n$,
\[
 \sum_m \sum_{j\ge m} (1-\big(1-2^{-2^{d(j-m)}} -\exp(-2^{-md} 2^{j\nu(\alpha_\text{min}+1/n)})\big)^{2^{md}}) <\infty.
\]
Therefore, the Borel-Cantelli lemma implies that for any $n$, there exists $m_0\in \N$ such that, for any $m\ge m_0$, any $r$ and any $j\ge m$, the inequality
\begin{equation}\label{eq:c-d}
 \| c^{(\cdot)}_{j,\cdot} -d^{(\cdot)}_{j,\cdot} \|_{\ell^\infty(T_{r,m})} \ge 2^{-j(\alpha_\text{min}+1/n)}
\end{equation}
holds almost surely. Using inequality~(\ref{eq:c-d}), we see that for any $n\in\N$,
\[
 \begin{array}{cc}
 \displaystyle \max\{
 \sup_{j\le l \le j+\log_2 j} \| c^{(\cdot)}_{l,\cdot} -d^{(\cdot)}_{l,\cdot} \|_{\ell^\infty(T_{r,m})}, &\\
 & \displaystyle \hspace{-100pt} 2^{-j([\alpha_\text{min}]+1)} \sup_{j-\log_2 j\le l\le j} (2^{l([\alpha_\text{min}]+1)} \|c^{(\cdot)}_{l,\cdot} -d^{(\cdot)}_{l,\cdot}\|_{\ell^\infty(T_{r,m})}
 \}
 \end{array}
\]
is larger than $2^{-j(\alpha_\text{min}+1/n)}$ for any $r,m$ almost surely. Moreover, since $X_\nu+f\in C^{\alpha_\text{min}}(\T^d)$ almost surely, we also almost surely have that, for any $j$ and any $r,m$,
\[
 \begin{array}{cc}
 \displaystyle \max\{
 \sup_{j\le l \le j+\log_2 j} \| c^{(\cdot)}_{l,\cdot} -d^{(\cdot)}_{l,\cdot} \|_{\ell^\infty(T_{r,m})}, \\
 &\displaystyle \hspace{-100pt} 2^{-j([\alpha_\text{min}]+1)} \sup_{j-\log_2 j\le l\le j} (2^{l([\alpha_\text{min}]+1)} \|c^{(\cdot)}_{l,\cdot} -d^{(\cdot)}_{l,\cdot} \|_{\ell^\infty(T_{r,m})}
 \}
 \end{array}
\]
is lower than $2^{-j\alpha_\text{min}}$. These two bounds allows us to say that for any $n\in\N$, the quantity
\[
 \lim_{j\to\infty}\frac{
 \begin{array}{cc}
 \displaystyle \log_2 \max\{
 \sup_{j\le l \le j+\log_2 j} \| c^{(\cdot)}_{l,\cdot} -d^{(\cdot)}_{l,\cdot} \|_{\ell^\infty(T_{r,m})}, \\
  & \displaystyle \hspace{-150pt} 2^{-j([\alpha_\text{min}]+1)} \sup_{j-\log_2 j\le l\le j} (2^{l([\alpha_\text{min}]+1)} \|c^{(\cdot)}_{l,\cdot} -d^{(\cdot)}_{l,\cdot} \|_{\ell^\infty(T_{r,m})}\}
 \end{array}
  }{-j}
\]
almost surely belongs to $[\alpha_\text{min},\alpha_\text{min}+1/n]$, for any $r,m$. Since this relation is valid for any $n$, Corollary~\ref{cor:carac} allows to conclude.
\EProof

\end{document}